\documentclass[a4paper,10pt,twoside]{amsart}
\usepackage[utf8]{inputenc}

\usepackage{amsfonts}
\usepackage{amsthm}
\usepackage{color}
\usepackage{graphicx}
\usepackage{amsmath,amssymb}

\def\R{\mathbb R}
\def\BB{\mathbb B} 
\def\SS{\mathbb S}

\def\Rm{\mathbb R^m}

\def\intO{\int_\Omega}

\def\intRm{\int_{\Rm}}
\def\intT{\int_0^T}

\def\eps{\varepsilon}

\def\esssup{\mathop{\rm ess\,sup\,}} % essential supremum
\def\diam{{\rm{diam}\,}}              % diameter of a set
\def\dist{{\rm{dist}\,}}              % distance
\def\div{{\rm{div}\,}}               % divergence of a vector field

\def\supp{{\rm{supp}\,}}              % support of a function

\newcommand{\norm}[1]{\left\|{#1}\right\|} %

\def\Xint#1{\mathchoice %  integral with something
{\XXint\displaystyle\textstyle{#1}}%
{\XXint\textstyle\scriptstyle{#1}}%
{\XXint\scriptstyle\scriptscriptstyle{#1}}%
{\XXint\scriptscriptstyle\scriptscriptstyle{#1}}%
\!\int}
\def\XXint#1#2#3{{\setbox0=\hbox{$#1{#2#3}{\int}$ }
\vcenter{\hbox{$#2#3$ }}\kern-.6\wd0}}

\def\dint{\Xint-} %mean integral

\newcommand{\pcz}[2]{\frac{\partial {#1}}{\partial {#2}}}
\newcommand{\pcza}[1]{\frac{\partial }{\partial {#1}}}

\newtheorem{theorem}{Theorem}[section]

\newtheorem{lemma}[theorem]{Lemma}

\theoremstyle{definition}

\newtheorem{remark}[theorem]{Remark}

\numberwithin{equation}{section}

\title[Conditional regularity for $p$-harmonic flows]{A conditional regularity result\\ for p-harmonic flows}\thanks{Revised version of \today}

\author[K. Kazaniecki]{Krystian Kazaniecki}
\address{Krystian Kazaniecki\newline \indent
Institute of Mathematics, University of Warsaw\newline \indent
Banacha 2, 02--097 Warszawa, Polska}
\email{K.Kazaniecki@mimuw.edu.pl}

\author[M. \L{}asica]{Micha\l{} \L{}asica}
\address{Micha\l{} \L{}asica\newline \indent
Institute of Applied Mathematics and Mechanics, University of Warsaw\newline \indent
Banacha 2, 02--097~Warszawa, Polska}
\email{M.Lasica@mimuw.edu.pl}

\author[K.E. Mazowiecka]{Katarzyna Ewa Mazowiecka}
\address{Katarzyna Ewa Mazowiecka\newline \indent
Institute of Mathematics, University of Warsaw\newline \indent
Banacha 2, 02--097 Warszawa, Polska}
\email{K.Mazowiecka@mimuw.edu.pl}

\author[P. Strzelecki]{Pawe\l\ Strzelecki}
\address{Pawe\l{} Strzelecki\newline \indent
Institute of Mathematics, University of Warsaw\newline \indent
Banacha 2, 02--097 Warszawa, Polska}
\email{P.Strzelecki@mimuw.edu.pl}

\date{}


\begin{document}
\maketitle

\begin{abstract}
We prove an $\varepsilon$-regularity result for a wide class of parabolic systems
$$
u_t-\text{div}\big(|\nabla u|^{p-2}\nabla u) = B(\,\cdot\,, u, \nabla u)
$$
with the right hand side $B$ growing critically, like $|\nabla u|^p$. It is assumed \emph{a priori} that the solution $u(t,\cdot)$ is uniformly small in the space of functions of bounded mean oscillation. The crucial tool is provided by a sharp nonlinear version of the Gagliardo--Nirenberg inequality which has been used earlier in the elliptic context by T.~Rivi\`{e}re and the last named author. 
\end{abstract}
% 
% \begin{figure}[h!]
% \centering
% \includegraphics[width=50mm]{Zeus-Releases-The-Kraken-6-4-10-kc}
% \caption{Liam Neeson i jego podejscie do teorii regularnosci }
% \label{fig:method}
% \end{figure}

\section{Introduction}

It is very well-known that for second order nonlinear parabolic and elliptic systems with the right hand side growing critically --- so that standard bootstrap methods do not yield any extra regularity of the solutions --- regularity is a subtle issue. For some systems (e.g. harmonic maps from planar domains into Riemannian manifolds or the classic $H$-surfaces) all the solutions are regular, for others --- classes of partially regular solutions, with some control of the size and structure of the singular set, coexist with wildly singular solutions. Again, harmonic maps into Riemannian manifolds provide a well-known example: minimizing and stationary solutions are partially regular, whereas on the other hand there are examples of everywhere discontinuous solutions; see e.g. \cite{hardt1997} and \cite{helein-book} for more. 

Regularity is often related to a special structure of the system (which allows one to exploit the compensation phenomena); a number of regularity results on harmonic or $p$-harmonic maps can be rephrased as follows: ``once the solution is small in the space $BMO$ of functions having bounded mean oscillation, it must be regular''.  In \cite{rivierestrzelecki}, T.~Rivi\`{e}re and the last named author have shown that $\eps$-regularity results of that type can be obtained \emph{without} any assumptions on the structure of a degenerate elliptic system. Of course, there is \emph{a small price to pay for beauty}\footnote{Familiarity with the topic, or with \emph{Butch Cassidy and the Sundance Kid}, reveals the catch.}: one must a priori know that the solution is small in $BMO$ \emph{and}  has second order distributional derivatives in $L^p$; on the other hand, even for the unconstrained $p$--Laplace system $\mathrm{div} (|\nabla u|^{p-2}\nabla u) = 0$ the $W^{2,p}$--regularity is not known.

In the present note, we provide a parabolic counterpart of the elliptic $\eps$-regularity result from \cite{rivierestrzelecki}. The main message is that if a solution $u(t,\cdot)$ is uniformly small in the space of functions of bounded mean oscillation \emph{and} has second order spatial derivatives, then it must be regular. Because of both assumptions, our main result is only \emph{conditional}, as the title indicates. (We hint below at a particular example where the first assumption alone implies the second one, but the whole situation is far from being fully understood.)

Let us now pass to more formal and precise statements.

We study the system of nonlinear parabolic equations,
\begin{equation} 
 u_t - \div (|\nabla u|^{p-2} \nabla u) = B(\,\cdot\, , u, \nabla u)\, ,
 \label{p-harmonic}
\end{equation}
for a vector $u = (u^1, \ldots u^N)$, given a vector $B=(B^1,\ldots,B^N)$ with 
\[
u^k \colon \mathbb (0,T] \times \Omega \to \mathbb R,\quad B^k \colon \Omega \times \R^{N} \times \R^{Nm} \to \R ,  \qquad 1\le k\le N,
\]
where $\Omega$ stands for an open domain in $\R^m$. We restrict ourselves to the case $p\ge 2$. We assume that the functions $B^k$, prescribing the nonlinearity of the right hand side, satisfy the critical growth condition
\begin{equation} 
 |B^k(\,\cdot\, ,u,\nabla u)| \leq \Lambda |\nabla u|^p . 
 \label{growth}
\end{equation}

The study of regularity of weak solutions to equations involving the parabolic $p$-Laplace operator $u_t - \div (|\nabla u|^{p-2} \nabla u)$ was initiated by DiBenedetto and Friedman \cite{dibenedetto2} (see also DiBenedetto's book \cite{dibenedetto3} and the recent memoirs of Duzaar, Mingione and Steffen \cite{duzaarmingionesteffen} and of B\"{o}gelein, Duzaar and Mingione \cite{bogeleinduzaarmingione}). Known explicit solutions to \eqref{p-harmonic} with $B^k = 0$ show that even in this case no regularity higher than $C^{1,\alpha}$ can be expected. Such regularity was indeed shown in \cite{dibenedetto2} for \eqref{p-harmonic} with $B^k$ satisfying at most
\begin{equation} 
 |B^k(\,\cdot\, ,u,\nabla u)| \leq \Lambda |\nabla u|^{p-1} . 
 \label{controllable}
\end{equation}

However, the more restrictive growth assumption \eqref{controllable} is typically \emph{not} satisfied in numerous inhomogeneous flows of geometric origin.  The $p$-harmonic heat flow with values in a Riemannian manifold $\mathcal{N}$ (which is assumed to be isometrically embedded in some ${\mathbb R}^d$) is a prototype of \eqref{p-harmonic}; in this case 
\begin{equation} 
 B(\,\cdot\, ,u,\nabla u) = |\nabla u|^{p - 2} A(u)(\nabla u, \nabla u), 
 \label{fundamental}
\end{equation}
where $A(u)(\cdot)$ is the second fundamental form of $\mathcal N$ at $u(\cdot)$. The interest in such systems arose in connection with the \emph{homotopy problem for ($p$-)harmonic maps}, i.\,e.\;the problem of finding a ($p$-)harmonic map homotopic to a given map between smooth compact Riemannian manifolds $\mathcal M$ and $\mathcal N$.\footnote{Strictly speaking, in this case the p-Laplace operator has to be substituted with appropriate p-Laplace-Beltrami operator on $\mathcal M$.} In this context, with $m=2,\ p=2$, the flow of harmonic maps was investigated by Eells and Sampson \cite{eellssampson}, who proved the existence of global regular solutions under the assumption that $\mathcal N$ has nonpositive sectional curvature (and solved the homotopy problem in this case). 

Without this assumption on the curvature of $\mathcal N$ one cannot expect global regularity. In the case $\mathcal N = \mathbb S^m,\ m>3$ examples of evolutions with finite time blow-up were constructed (see \cite{coronghidaglia}, also \cite{changdingye} and 
references therein). Later, Struwe \cite{struwe1,struwe2} and Struwe in collaboration with Chen \cite{chenstruwe} considered weak solutions to the harmonic heat flow with values in arbitrary $\mathcal N$ and were able to show existence and bounds on the set of singularities. The crucial tool in their work (in case $m>2$) was Struwe's monotonicity formula \cite[Proposition 3.3.]{struwe2}, which is only available if $p = 2$. 

To the best of our knowledge, there is no fully general existence result for the flow of $p$-harmonic maps into an arbitrary compact manifold for $p > 2$. Partial results in this direction include the  work of Misawa \cite{misawa1}
%and Fardoun and Regbaoui \cite{fardounregbaoui}
which generalizes the classic work of Eells and Sampson ($p=2$) obtaining global regular solutions for $\mathcal N$ of non-positive sectional curvature, and the paper \cite{hunger2} by Hungerb\" uhler, where the existence of weak solutions is shown for $\mathcal N$ being a homogeneous space. In the latter case, the proof exploits symmetry of the image and no additional regularity is obtained. Finally,  Hungerb\" uhler \cite{hunger1}  obtained the existence of global weak solutions (regular except a finite set of time instances) in the conformally invariant case $p = m$. (Also for $p=m$; B\"{o}gelein, Duzaar and Scheven \cite{bogeleinduzaarscheven} extend the results of \cite{hunger1} to domains with boundary, and construct a solution with H\"{o}lder continuous spatial gradient in the case when the target $\mathcal N$ has non-positive sectional curvature.) The proof in \cite{hunger1} is based on a local estimate on energy concentration together with a conditional a priori estimate that controls the norms of $\nabla u$ in higher $L^q$ spaces allowing to bound $\nabla u$ in $L^\infty$ via Moser iteration. Those estimates hold provided the $m$-energy is appropriately small,
\begin{equation}\label{hungcond}
 \sup_{t,x} \int_{B_R(x)} |\nabla u(t,x)|^m dx <\varepsilon.
\end{equation}          

In the present paper, we obtain conditional estimates of a similar form, but instead of smallness of the energy \eqref{hungcond}, we assume smallness of local $BMO$ norm of the solution (which, in light of the embedding $W^{1,m}$ into $BMO$ in dimension $m$, is a weaker assumption). In particular, our method works for any $p>2$ and (formally) we only need to control a norm of the solution and not of its derivatives. On the other hand, we have no proof of existence; therefore we do need the assumption that the solution actually exists in $C((0,T], L^2(\Omega)) \cap L^p((0,T],W^{2,p}(\Omega))$. Let us now state our results. 

\begin{theorem} 
 Assume that $u \in C((0,T], L^2(\Omega)) \cap L^p((0,T], W^{2,p}(\Omega))$, where $\Omega$ is an open domain in $\R^m$, is a weak solution to \eqref{p-harmonic}. Let $q > p$, $\Omega' \subset \subset \Omega$, $\delta > 0$. There exist a positive number $\varepsilon_0 = \varepsilon_0(\Omega,N,p,q,\Lambda)$ such that if 
\begin{equation}
	\label{bmo-small}
	\|u\|_{L^\infty((0, T],BMO(\Omega))} < \varepsilon_0
\end{equation} 
 then there holds 
 $$\|\nabla u \|_{L^q((\delta, T]\times \Omega')} \leq C,$$
 where $C=C(\Omega, N, p, q, \Lambda, \dist(\Omega', \partial \Omega), \delta)$.
 \label{lq}
\end{theorem}    

Together with the work of DiBenedetto and Friedman, this theorem implies 
the following. 
\begin{theorem}    Assume that $u \in C((0,T], L^2(\Omega)) \cap L^p((0,T], W^{2,p}(\Omega))$ is a weak solution to \eqref{p-harmonic}. There exists a constant $\varepsilon_0 = \varepsilon_0(\Omega,m,p,\Lambda)$ such that if the condition 
 $$\|u\|_{L^\infty((0,T],BMO(\Omega))} < \varepsilon_0$$ 
 is satisfied, then $u \in C^{1, \alpha}_\mathrm{loc}((0,T] \times \Omega)$. 
 \label{regularity}
\end{theorem}

\begin{remark}
\label{duzremark}	
	 (1) It would be desirable to prove the above results under a somewhat different assumption that one deals with a weak solution
	\[
	u \in L^p((0,T], W^{1,p}(\Omega)) \qquad\mbox{such that}\qquad 
	|\nabla u|^{(p-2)/2} \nabla u \in L^2((0,T], W^{1,2}(\Omega))\, 
	\]
which satisfies the smallness condition \eqref{bmo-small}. 
One reason for this is that even in the case of $B=0$, i.e. in the case of the pure parabolic $p$-Laplacian the existence of second spatial derivatives in $L^p$ is unknown, whereas the existence of first spatial order derivatives of $|\nabla u|^{(p-2)/2} \nabla u$ can be shown, again for $B=0$, by the standard method of difference quotients.\footnote{For this remark, we are indebted to Frank Duzaar and a discussion with him during the meeting of German and Polish Mathematical Societies.}

However, despite some efforts, we were not able to achieve this. One of the technical difficulties behind the failure is hidden in the proof of the Caccioppoli inequality.

(2) On the other hand, we wish to observe that in some special cases, e.g. when $n=p=2$ and $B\not=0$ has a div--curl structure, as in the case of the flow of $H$-surfaces with constant $H$, 
\[
u_t = \Delta u -2H u_x\wedge u_y,
\]
where $2Hu_x\wedge u_y$ is a Jacobian-like term, the smallness condition \eqref{bmo-small} \emph{alone} is sufficient to prove that $\nabla u$ is in $L^2$--$W^{1,2}$ locally (to this end, one employs the familiar methods based on the duality of Hardy space and $BMO$). In other words, in these cases one could remove the assumption on $D^2 u$ from the analogues of both of our theorems. 

We do not know if there are any other, significantly different, examples of a similar phenomenon where smallness of a solution in $L^\infty$--$BMO$ \emph{implies} the existence of second order spatial derivatives of this solution.
\end{remark}

\smallskip

There are other conditional results, involving smallness assumptions in \emph{stronger} norms but with no a priori assumption on the second derivatives of the solution. In this direction, Leone, Misawa and Verde \cite{misawa4} give a natural smallness condition in the $L^\infty$--norm, yielding an estimate of the parabolic Hausdorff dimension of the singular set. Karim and Misawa \cite{misawa5} consider the case of $\frac{2m}{m+2}<p<2$.

Finally, comparing our results with Kuusi and Mingione \cite{kuusimingione}, we see that Corollary~\ref{regularity} is, on one hand, stronger than e.g.\;\cite[Thm~1.3]{kuusimingione} which gives a borderline version of the $L^\infty$--boundedness of $Du$ from \cite[Chapter VIII]{dibenedetto3}. Even for solutions $u\in L^\infty((0,T], BMO\cap W^{2,p}(\Omega))$ the right hand side of \eqref{p-harmonic} is formally only in $L^2$, and not in the Lorentz space $L^{m+2,1}$, as it is assumed in \cite[Thm~1.3]{kuusimingione}. On the other hand, the smallness assumption \eqref{bmo-small} is pretty strong and seems to be restrictive -- yet, in a sense, necessary: even in the elliptic case, e.g.\;for harmonic maps into Riemannian manifolds, an assumption of this kind is needed for regularity; it is well known that near an isolated singularity of a harmonic map $w\colon \BB^3\to\SS^2$ the $BMO$ norm of $w$ is \emph{not} small. 

Our main technical tool is an interpolation inequality of Gagliardo-Nirenberg type, discovered by Rivi\`{e}re and the last named author of the present paper, see \cite{rivierestrzelecki}. Let $\psi \in C^\infty_c(\mathbb R^m)$ be fixed. Using $\mathcal H^1 - BMO$ duality, the authors of \cite{rivierestrzelecki} proved the existence of a constant $C=C(m)$ such that
\begin{equation}
\begin{aligned}
\intRm \psi^{s+2} |\nabla u|^{s+2}
& \le  C s^2 \|u\|_{BMO(\Rm)}^2 \biggl\{
\intRm \psi^{s+2} |\nabla u|^{s-2}
|\nabla^2 u|^2
\biggr. \\
&\quad {}\biggl.
\phantom{C s^2 \|u\|_{BMO(\Rm)}^2}
+ \|\nabla\psi\|_{L^\infty(\Rm)}^2
\intRm \psi^{s} |\nabla u|^{s}\biggr\}
\end{aligned} 
\label{int-bmo}
\end{equation}
for any function $u\in W^{2,1}_{{\rm loc}}(\Rm) \cap BMO(\Rm)$ for which the right hand side is finite. A version of \eqref{int-bmo} in time-dependent setting follows immediately. Let now $\psi \in C^\infty_c((0,T]\times\Rm)$. Integrating \eqref{int-bmo} over time yields 
\begin{equation}
\begin{aligned}
&\intT\!\! \intRm \psi^{s+2} |\nabla u|^{s+2} \le  C s^2 \|u\|_{L^\infty((0,T],BMO(\Rm))}^2 
\\ &\quad \cdot\biggl\{\intT \!\!\intRm \psi^{s+2} |\nabla u|^{s-2}
|\nabla^2 u|^2 + \|\nabla\psi\|_{L^\infty((0,T]\times\Rm)}^2
\intT \!\!\intRm \psi^{s} |\nabla u|^{s}\biggr\}
\end{aligned} 
\label{int-bmot}
\end{equation}
for any $u\in L^1((0,T], W^{2,1}_{{\rm loc}}(\Rm)) \cap L^\infty((0,T], BMO(\Rm))$ such that the right hand side is finite. This inequality allows us to control the right hand side of \eqref{p-harmonic} provided that the smallness condition from Theorem \ref{lq} is satisfied. This part of our work is in fact a parabolic version of \cite{rivierestrzelecki}. Having obtained such conditional bounds on sufficiently high local $L^q$ norms of the gradient, one may then use known results of \cite{dibenedetto3} to obtain Theorem~\ref{regularity} (see Section~3.2). 

\bigskip
\noindent
{\bf Notation.} We denote by $BMO(\Omega)$ the space of functions on a given domain $\Omega$ of \emph{bounded mean oscillation}, with the seminorm
\[
\norm{f}_{BMO(\Omega)}:=\sup_Q\left(\dint_Q|f(y)-f_Q|dy\right)<\infty,
\]
the supremum being taken over all cubes $Q$ in $\Omega$, where $f_Q
$ denotes the average of $f$ on $Q$, i.e. $\dint_Q f dx = |Q|^{-1}\int_Q f dx$, $|Q|$ being the Lebesgue measure of $Q$. Given $t \in (0,T]$, we write and $\Omega_t$ for the cylindrical domain $(0,t]\times\Omega$.

\section{Caccioppoli inequality}

In order to obtain $L^q$ estimates in the homogeneous case (i.\,e.\;the right hand side of \eqref{p-harmonic} equal $0$) DiBenedetto and Friedman \cite{dibenedetto1} tested the system with $\zeta^2|\nabla u|^{2\alpha}\nabla u$, where $\zeta$ is a suitably chosen smooth cutoff function. The same standard test functions were used by Hungerb\"{u}hler in the conformally invariant case and by Rivi\`{e}r{e} and the last author in the elliptic case. 

We modify their derivation obtaining the following Caccioppoli inequality for derivatives of solutions of (\ref{p-harmonic}).

\begin{lemma} Assume that $u\in C((0,T], L^2(\Omega)) \cap L^p((0,T],W^{2,p}(\Omega))$ is a weak solution of
(\ref{p-harmonic}). Let $\zeta\in C_c^\infty((0,T]\times \Omega)$ and set
$w\colon\!\!\! = |\nabla u|^2$. There exists a constant
$C_1=C_1(m,N,p,\Lambda)$ such that for each $\alpha\ge 0$ we have
\begin{equation}
\begin{split}
\frac{1}{2+ 2 \alpha} & \esssup_{t \in (0,T]} \intO \zeta(t, \cdot)^2 w^{1+ \alpha}  + 
\frac{p-2+\alpha}{8} \iint_{\Omega_T}  \zeta^2
w^{\frac{p}{2} -2 + \alpha} \left|\nabla w\right|^2 \\
 + &\frac 12\iint_{\Omega_T} \zeta^2 w^{\frac p2 -1+\alpha} |\nabla^2
u|^2 \le \biggl(\frac{2(p-1)^2}{p-2+2\alpha}+\frac 12\biggr)
\iint_{\Omega_T} |\nabla \zeta|^2
w^{\frac{p}{2}+\alpha} \\
&+ \frac{1}{1 + \alpha} \iint_{\Omega_T} \zeta \zeta_t w^{1+\alpha} + C_1 (p+\alpha) \iint_{\Omega_T} \zeta^2 w^{\frac p2 +1+\alpha},
\end{split}
\label{caccio}
\end{equation}
provided the right hand side is finite.
\end{lemma}

\noindent
{\bf Proof.}
We would like to take the spatial derivative of \eqref{p-harmonic} and test the resulting equation with $\zeta^2|\nabla u|^{2\alpha} \nabla u$. However, this is not an admissible test function due to lack of time regularity. To circumvent this difficulty, we exploit the standard technique of time mollification by Steklov averages. 

First, note that there exists $\delta \in (0,T)$ such that $\supp \zeta \in (2\delta,T]\times\Omega$. For any given $0< h < \delta$, by $I_h$ we denote the operator of Steklov average, i.\,e.
$$I_h(u)(t,\cdot) = \frac{1}{h}\int_{t - h}^t u.$$ 
If $X$ is a Banach space, $q> 1$, $t>0$, $I_h$ maps $L^q((0,t], X)$ to $W^{1,q}((\delta,t], X)$ and $I_h(u) \to u$ in $L^q((\delta,t], X)$ as $h \to 0^+$. 

We apply $I_h$ to the weak formulation of (\ref{p-harmonic}) and differentiate both sides with respect to $x_j$. Testing the resulting equation with $\varphi_h^{ij} = \zeta^2|I_h \nabla u|^{2\alpha} I_h \pcz{u^i}{x_j}$ we obtain, for any $t>\delta$ (indices $i,j$ are summed),

\begin{multline}
\label{steklov} \iint_{\Omega_t} \pcza{t} \left(\pcz{I_h u^i}{x_j}\right) \varphi_h^{ij} \\ {}+
\iint_{\Omega_t} I_h\left[ |\nabla u|^{p-2}
\nabla\left(\pcz{u^i}{x_j}\right) + \pcza{x_j}\left(|\nabla u
|^{p-2}\right)\nabla u^i\right]\cdot \nabla\varphi_h^{ij} \\
=-\iint_{\Omega_t} \pcz{\varphi_h^{ij}}{x_j} I_h B^i(x,u,\nabla u) .
\end{multline}

We first deal with the term in \eqref{steklov} containing the time derivative. We have 
\begin{multline} 
\label{steklov2} 
 \iint_{\Omega_t} \pcza{t} \left(\pcz{I_h u^i}{x_j}\right) \varphi_h^{ij} = \frac{1}{2 + 2 \alpha} \iint_{\Omega_t}  \pcza{t} \left(|I_h \nabla u|^{2 + 2 \alpha}\right) \zeta^2 \\= \frac{1}{2 + 2 \alpha} \int_\Omega  |I_h \nabla u|^{2 + 2 \alpha} \zeta^2(t, \cdot) - \frac{1}{1 + \alpha} \iint_{\Omega_t}  |I_h \nabla u|^{2 + 2 \alpha} \zeta \pcz{\zeta}{t} .
\end{multline}

The second term on the right hand side of \eqref{steklov2} clearly converges as $h\to 0^+$ to respective non-averaged expression, as do the last two integrals in \eqref{steklov}. As all of them are bounded independently of $t \in (\delta, T]$, we deduce from \eqref{steklov} and \eqref{steklov2} that $t\mapsto \int_\Omega  |\nabla u|^{2 + 2 \alpha} \zeta^2(t,\cdot)$ belongs to $L^\infty((\delta, T])$ and,  for almost every $t \in (\delta, T]$, is equal to the limit (as $h\to 0$) of the first integral on the right hand side of \eqref{steklov2}. Eventually, passing to the limit in \eqref{steklov}, we obtain 
for almost all $t\in (\delta,T]$ 
\begin{multline}
\label{twotwo} \frac{1}{2 + 2 \alpha} \int_\Omega  |\nabla u|^{2 + 2 \alpha} \zeta^2(t, \cdot) - \frac{1}{1 + \alpha} \iint_{\Omega_t}  |\nabla u|^{2 + 2 \alpha} \zeta \pcz{\zeta}{t} \\ {}+
\iint_{\Omega_t} \left[ |\nabla u|^{p-2}
\nabla\left(\pcz{u^i}{x_j}\right) + \pcza{x_j}\left(|\nabla u
|^{p-2}\right)\nabla u^i\right]\cdot \nabla\varphi^{ij} \\
=-\iint_{\Omega_t} \pcz{\varphi^{ij}}{x_j} B^i(x,u,\nabla u), 
\end{multline}
where $\varphi^{ij}$ denotes $\zeta^2|\nabla u|^{2\alpha} \pcz{u^i}{x_j}$.

We now estimate the left and right hand side of \eqref{twotwo} separately.

\smallskip
\noindent {\em Left hand side of {\em (\ref{twotwo}).}}
 A~routine but somewhat tedious computation leads to the following two
equalities: 
\begin{equation}
\begin{aligned}
\iint_{\Omega_t} |\nabla u|^{p-2} \nabla\left(\pcz{u^i}{x_j}\right)
\cdot  \nabla \varphi^{ij} = \iint_{\Omega_t} \zeta^2
w^{\frac{p-2}{2}+\alpha}
\left|\nabla\left(\pcz{u^i}{x_j}\right)\right|^2\\
+ \frac{\alpha}{2} \iint_{\Omega_t} \zeta^2
 w^{\frac{p-2}{2}-1+\alpha} |\nabla w|^2
 + \iint_{\Omega_t} \zeta(\nabla\zeta\cdot\nabla w)
 w^{\frac{p-2}{2}+\alpha}\\
 =: I_1+I_2+I_3\, ; 
\end{aligned}
\label{lhs-1}
\end{equation}
\begin{equation}
\begin{aligned}
\iint_{\Omega_t} \pcza{x_j}& \left(|\nabla u |^{p-2}\right)\nabla
u^i\cdot \nabla\varphi^{ij}\\
& = \frac{p-2}{4} \iint_{\Omega_t}
\zeta^2 w^{\frac{p-2}{2}-1+\alpha} |\nabla w|^2\\
&\quad{}+ \frac{(p-2)\alpha}{2} \iint_{\Omega_t}
\zeta^2w^{\frac{p-2}{2}-2+\alpha} \sum_i\left(\nabla w\cdot \nabla
u^i\right)^2\\
&\quad {}+ (p-2) \iint_{\Omega_t}\zeta w^{\frac{p-2}{2}-1+\alpha}
\sum_i (\nabla\zeta\cdot \nabla u^i)(\nabla w\cdot \nabla u^i)\\
& =: I_4+I_5+I_6 .
\end{aligned}
\label{lhs-2}
\end{equation}
Using Cauchy's inequality with $\varepsilon$: $ab\le
\frac{\eps^2a^2}{2}+\frac{b^2}{2\eps^2}$, we estimate
\begin{equation}
\begin{split}
|I_3|+|I_6| & \le  (p-1) \iint_{\Omega_t} \zeta|\nabla\zeta|\, |\nabla
w| w^{\frac{p-2}{2}+\alpha} \\
& \le \frac{(p-1)\eps^2}{2} \iint_{\Omega_t} \zeta^2
w^{\frac{p-2}{2}-1+\alpha}|\nabla w|^2 +\frac{p-1}{2\eps^2}
\iint_{\Omega_t} |\nabla \zeta|^2 w^{\frac{p}{2}+\alpha}\, ,
\end{split}
\label{absorb}
\end{equation}
splitting the integrand so that the term with $|\nabla w|^2$ can be absorbed in $I_2+I_4$. Choosing $\eps^2$ so that $(p-1)\eps^2/2 = (p-2+2\alpha)/8$, and
combining (\ref{lhs-1}), (\ref{lhs-2}) and (\ref{absorb}), we
obtain finally
\begin{equation}
\begin{split}
\mbox{left hand side of (\ref{twotwo})} &\ge \frac{1}{2+ 2 \alpha} \intO w^{1+ \alpha} \zeta^2(t,\cdot) - \frac{1}{1+ \alpha} \iint_{\Omega_t} w^{1+ \alpha} \zeta \pcz{\zeta}{t} \\ 
&\quad{}+\frac{p-2+2\alpha}8
\iint_{\Omega_t} \zeta^2 w^{\frac{p-2}{2}-1+\alpha}|\nabla w|^2 \\
&\quad{}+ \iint_{\Omega_t} \zeta^2w^{\frac{p-2}{2}+\alpha} |\nabla
u^i_{x_j}|^2 \\
&\quad{}+\frac{(p-2)\alpha}{2} \iint_{\Omega_t}
\zeta^2w^{\frac{p-2}{2}-2+\alpha} \sum_i\left(\nabla w\cdot \nabla
u^i\right)^2 \\
&\quad{}-\frac{2(p-1)^2}{p-2+2\alpha} \iint_{\Omega_t} |\nabla \zeta|^2
w^{\frac{p}{2}+\alpha} .
\end{split}
\label{lhs-final}
\end{equation}

\smallskip
\noindent {\em Right hand side of  {\em(\ref{twotwo})}.}
Using the growth condition $|B(x,u,\nabla u)|\le \Lambda |\nabla
u|^p$, we write
\begin{equation}
\biggl|\iint_{\Omega_t} \pcz{\varphi^{ij}}{x_j} B^i(x,u,\nabla u)
\biggr|\le C (J_1+J_2+J_3)\, ,
\end{equation}
where the constant $C=C(m,N,\Lambda)$ and
\begin{equation}
\begin{split}
J_1 = \iint_{\Omega_t} \zeta^2 w^{\frac p2 +\alpha} \bigl|\nabla
u^i_{x_j}\bigr|  \qquad
J_2 = \alpha \iint_{\Omega_t} \zeta^2 w^{\frac{p-1}2 +\alpha} \left|\nabla
w\right| \\
J_3 = \iint_{\Omega_t} \zeta|\nabla\zeta| w^{\frac{p+1}2 +\alpha} .
\end{split}
\end{equation}
Set
\[
J_0\colon= \iint_{\Omega_t} \zeta^2w^{\frac p2 + 1 +\alpha}\, .
\]
To absorb all terms that contain second order derivatives of $u$,
we again apply the Cauchy--Schwarz inequality in
a familiar way and obtain
\begin{eqnarray}
\nonumber
J_1 & \le & \frac{\eps_1^2}{2} \iint_{\Omega_t} \zeta^2
w^{\frac{p-2}{2} +\alpha} \bigl|\nabla u^i_{x_j}\bigr|^2 
 + \frac{1}{2\eps_1^2} J_0\, ,\\
\nonumber
J_2 & \le & \frac{\alpha\eps_2^2}{2} \iint_{\Omega_t} \zeta^2
w^{\frac{p}{2} -2+\alpha} \left|\nabla w\right|^2 
+ \frac{\alpha}{2\eps_2^2} J_0\, .
\end{eqnarray}
Finally,
\begin{equation}
J_3\le \frac C2 J_0 + \frac 1{2C} \iint_{\Omega_t} |\nabla\zeta|^2
w^{\frac p2 + \alpha}\, .
\end{equation}
Making appropriate choices of $\eps_1,\eps_2>0$, we combine
the estimates of $J_1,J_2,J_3$ with (\ref{lhs-final}) and, taking supremum of both sides over $t\in (0,T]$, complete the proof of the lemma.\hfill $\Box$

\section{Gradient estimates}

This section contains proofs of Theorems \ref{lq} and \ref{regularity}. The high integrability of solution to \eqref{p-harmonic} that we seek will eventually follow upon the iterations of the Caccioppoli inequality \eqref{caccio} combined with a parabolic version of the Sobolev inequality, 
 %\begin{equation}
 \begin{align}
   \nonumber    
  \iint_{(\delta,T]\times\Omega'} |w|^{\frac{p + 2 \alpha}{2} + \frac{2}{m}(1+\alpha)} &\leq C(m)\esssup_{t \in (\delta,T]}   \left(\int_{\Omega'} |w(t,\cdot)|^{1+\alpha}\right)^{\frac{2}{m}}\\\label{emb}
  &\quad \cdot\left(\iint_{(\delta,T]\times\Omega'} \big|\nabla w^{\frac{p + 2 \alpha}{4}}\big|^2
+ c_{\Omega'}\iint_{(\delta,T]\times\Omega'}|w|   ^{\frac{p+2\alpha}{2}}\right),
  \end{align}
 %\end{equation} 
 which holds for any $\delta > 0$, any smooth bounded $\Omega' \Subset \Omega$, and will be applied to $w=|w| = |\nabla u|^2$. Recall that inequality \eqref{emb} is obtained by an application of H\" older inequality and Sobolev embedding\footnote{We tacitly assume here that $m>2$; for $m=2$ the exponent $2m/(m-2)$ can be replaced by any $s\in (2,\infty)$.} $W^{1,2} \hookrightarrow L^{\frac{2m}{m-2}}$ in dimension $m$,  
\begin{align*} 
 \int_{\Omega'} |w|^{\frac{p + 2 \alpha}{2} + \frac{2}{m}(1 + \alpha)} &\leq \left(\int_{\Omega'} |w|^{1 + \alpha}\right)^{\frac{2}{m}} \left(\int_{\Omega'} |w|^{\frac{p + 2 \alpha}{4} \cdot \frac{2m}{m-2}} \right)^{\frac{m-2}{2m} \cdot 2} \\ &\leq  C(m) \left(\int_{\Omega'} |w|^{1 + \alpha}\right)^{\frac{2}{m}} \left(\int_{\Omega'} \big|\nabla w^{\frac{p + 2 \alpha}{4}}\big|^2+
c_{\Omega'}\int_{\Omega'} |w|^{\frac{p+2\alpha}{2}}\right) ,
\end{align*}
and integrating the result over the time interval $(\delta,T]$.  We also note that if $\Omega'$ is a cube or a ball, then the above inequalities hold with
\begin{equation}
	\label{comegaprim} 
	c_{\Omega'}=(\diam \Omega')^{-2}\, .
\end{equation}
 
 However, it is easy to see that if the exponent $\alpha$ is small, then the combined inequalities (Caccioppoli and `parabolic' Sobolev) do not yield any increase in integrability of $w$. This is caused precisely by the critical term \eqref{badterm}, with $|\nabla u|^{p+2+ 2\alpha}=w^{\frac p2 + 1+ \alpha}$. Therefore, we need first to control $\nabla u$ in a sufficiently high $L^q$ space. This is why some smallness assumption is required for our assertions to hold. We exploit such an assumption through a bootstrap procedure involving a finite number of applications of the interpolation inequality \eqref{int-bmot} to control the `bad' term 
 \begin{equation} 
 C_1(p+\alpha)\iint_{\Omega_T} \zeta^2 w^{\frac p2 +1+\alpha} 
 \label{badterm}
 \end{equation}
 on the right hand side of the Caccioppoli inequality. 
 
Once the threshold level of $q = p(m+2)$ is exceeded, the assertion of Theorem \ref{regularity} may be deduced from DiBenedetto's book. 

We set forth all the details below, starting with a~parabolic version of the `elliptic' argument from \cite{rivierestrzelecki}.

\subsection{$L^q$ estimate for $\nabla u$}

We now explain how to iterate the Caccioppoli inequality, using the Gagliardo--Nirenberg inequality at each step. Let $C_1$ and $C_2$ denote the constants from Caccioppoli inequality \eqref{caccio} and the interpolation inequality \eqref{int-bmot}, respectively. Fix  a sufficiently large number $\alpha_{{\text{max}}}$ that shall be specified later. 

We need the following smallness condition: 
 \begin{equation}
  C_1 C_2( p+2\alpha)^3 \|u\|_{L^{\infty}((0, T], BMO(\Omega))}^2 \leq \frac{1}{2},
 \label{small}
 \end{equation}
 for every $\alpha\in [0,\alpha_{{\text{max}}}]$. We choose two nonnegative functions $\zeta, \psi \in C^{\infty}_{c}(\Omega_T)$ so that $\psi\equiv 1$ on $(\delta,T]\times\Omega'$ and
 \begin{equation}\label{testujaca}
  \zeta^2=\psi^{p+2+2\alpha},\quad 0\le\zeta\le1, \quad |\nabla \zeta|\leq 2 \dist(\Omega', \partial \Omega)^{-1}, \quad  |\zeta_t|\leq 2 \delta^{-1}.
 \end{equation}
 % $\zeta^2=\psi^{p+2\alpha+2}$, $\psi\equiv 1$ on $(\delta,T]\times\Omega'$, $|\nabla \zeta|\leq 2 \dist(\Omega', \partial \Omega)^{-1}$ and $|\zeta_t|\leq 2 \delta^{-1}$.
 We set $s=p+2\alpha$ in the interpolation inequality \eqref{int-bmot} and use it to estimate the bad term \eqref{badterm}. 

Due to smallness condition \eqref{small} we get 
\begin{align} 
\frac{1}{2 + 2 \alpha}  \esssup_{t \in (0,T]} &  \int_\Omega \zeta(t, \cdot)^2 w^{1+ \alpha}  + \frac 12\iint_{\Omega_T} \zeta^2 w^{\frac p2 -1+\alpha} |\nabla^2 u|^2    \nonumber
\\  & \quad {}+ \frac{p-2+\alpha}{8} \iint_{\Omega_T}  \zeta^2 w^{\frac{p}{2} -2+\alpha} \left|\nabla w\right|^2 \nonumber    
\\ 
\label{ready-to-absorb}
& \le  \biggl(\frac{2(p-1)^2}{p-2+2\alpha}+\frac 12\biggr) \iint_{\Omega_T}|\nabla \zeta|^2 w^{\frac{p}{2}+\alpha} + \frac{1}{1+\alpha} \iint_{\Omega_T} \zeta \zeta_t w^{1+\alpha}
\\&  \nonumber{}\quad  +\frac{1}{2} \|\nabla\psi\|_{L^\infty(\Omega_T)}^2 \iint_{\Omega_T} \psi^{p+2\alpha} w^{\frac{p}{2}+\alpha} +\frac{1}{2}\iint_{\Omega_T} \zeta^2 w^{\frac{p}{2}-1+\alpha}|\nabla^2 u|^2.
\end{align}
The last term in the right-hand side cancels with the matching term in the left hand side. Invoking the properties of $\zeta$ and $\psi$, cf. \eqref{testujaca}, and remembering that 
$$\left(\frac{p+2\alpha}4\right)^2 w^{\frac{p}{2} -2+\alpha} \left|\nabla w\right|^2 = \left|\nabla w^{\frac{p+2\alpha}{4}}\right|^2,$$ 
we obtain from~\eqref{ready-to-absorb}, upon multiplication by $2(1+\alpha)$,                                         
\begin{equation}
	\label{estm2}     
 \begin{aligned}
 \esssup_{t \in (\delta,T]} \int_{\Omega'} w^{1+\alpha} &+ \iint_{(\delta,T]\times\Omega'} \left|\nabla w^{\frac{p+2\alpha}{4}}\right|^2
\\ &\le\gamma (1+\alpha) \biggl( \iint_{\Omega_T} w^{\frac{p}{2}+\alpha} +  \iint_{\Omega_T} w^{1+\alpha}\biggr),
\end{aligned}
\end{equation}
where $\gamma$ stands for a generic constant which may depend on $\Omega$, $T$, $m$, $p$, $\delta$, and $\Omega'$. Inequality  \eqref{estm2}, combined with the parabolic Sobolev embedding
\eqref{emb}, leads to the  estimate
\begin{equation} 
\begin{aligned}
 \biggl(\iint_{(\delta,T]\times\Omega'} |\nabla u|&^{p + 2 \alpha + \frac{4}{m}(1 + \alpha)}\biggr)^{\frac{1}{\kappa}}\\ &\leq 
 \gamma(1+\alpha)\biggl(
 \iint_{\Omega_T}
|\nabla u|^{p+2\alpha} +
 \iint_{\Omega_T}|\nabla u|^{2+2\alpha}\biggr),
 \end{aligned}
\end{equation}
where $\kappa=1+\frac{2}{m}$ (and $\gamma$ could change). Now, applying H\"older inequality to replace $2+2\alpha$ by $p+2\alpha$ in the second exponent of the right hand side, we obtain, adjusting the constant $\gamma$ again,
\begin{equation}
\begin{aligned}
 \bigg(\iint_{(\delta,T]\times\Omega'} & |\nabla u|^{p + 2 \alpha + \frac{4}{m}(1+\alpha)}\bigg)^{\frac{1}{\kappa}}\\ & \leq
 \gamma(1+\alpha) \biggl( \iint_{\Omega_T} |\nabla u|^{p+2\alpha} +
 \left(\iint_{\Omega_T}|\nabla u|^{p+2\alpha}\right)^{\frac{2+2\alpha}{p+2\alpha}}\biggr)\, .
\end{aligned}
\end{equation}
We iterate this inequality finitely many times, starting from $\alpha=0$. After each step we obtain higher local integrability of $\nabla u$ (on a smaller domain).  To achieve the desired goal, we choose $\alpha_{{\text{max}}}$ large enough to get $\nabla u\in L^q_{{\text{loc}}}(\Omega_T)$ for the given value of $q$. Since the time instant $\delta>0$ above can be arbitrary, Theorem \ref{lq} follows.

\subsection{$L^\infty$ estimates for $\nabla u$}   
By Theorem \ref{lq}, we have $|\nabla u|\in L^{q'}$ for some $q'>p(m+2)$; 
thus,  setting
 $C_1=0$ and $\varphi_2:=\Lambda |\nabla u|^p$ in the structure condition $(S_5)$ in \cite[Chapter VIII, p. 217--218]{dibenedetto3}, we obtain $\varphi_2^2\in L^q$ for $q=q'/(2p)>(m+2)/2$, so that the structure assumption $(S_6)$ in \cite[Chapter VIII, p. 217--218]{dibenedetto3} is satisfied. Now,  \cite[Lemma 4.2, p.234]{dibenedetto3} implies that $\nabla u$ is locally 
bounded, and by  \cite[Theorem 1.1', p.256]{dibenedetto3} $\nabla u$ is locally H\"older continuous. Theorem~1.2 follows.

\subsection*{Acknowledgements.} The authors are grateful to Frank Duzaar for his comments on an earlier version of this paper (see Remark~\ref{duzremark}). They are also grateful to the anonymous referee, for his/her remarks that helped to improve the presentation and fix a gap in the proof of the Caccioppoli inequality.

The first three named authors would like to thank Piotr B. Mucha for acquainting them with the subject of regularity for parabolic $p$-Laplace equation. 

The work of K. Mazowiecka and P. Strzelecki on this project has been partially supported by the NCN grant no. 2012/07/B/ST1/03366.

\nocite{misawa4}
\nocite{misawa5}
\nocite{bogeleinduzaarscheven}
\small
%\subsection{H\"older continuity of gradient}
%\input{holder}
\bibliographystyle{plain}
\bibliography{DeGiorgi}

\end{document}